\documentclass[thmsa,11pt]{article}
\usepackage{amssymb}

\input{tcilatex}
\newtheorem{Theorem}{Theorem}[section]
\newtheorem{Proposition}[Theorem]{Proposition}
\newtheorem{Corollary}[Theorem]{Corollary}
\newtheorem{Remark}[Theorem]{Remark}
\newtheorem{Lemma}[Theorem]{Lemma}
\newtheorem{Definition}[Theorem]{Definition}
\newtheorem{Construction}[Theorem]{Construction}

\begin{document}

\author{MARIA JOI\c{T}A}
\title{INDUCED\ REPRESENTATIONS\ OF\ LOCALLY\ \ $C^{*}$ -ALGEBRAS}
\maketitle

\begin{abstract}
In this paper, by analogy with the case of $C^{*}$-algebras, we define the
notion of induced representation of a locally $C^{*}$-algebra and then we
prove the imprimitivity theorem for induced representations of locally $%
C^{*} $-algebras.

MSC: 22D30; 46L05; 46L08
\end{abstract}

\section{Introduction}

Locally $C^{*}$-algebras generalize the notion of $C^{*}$-algebra. A locally 
$C^{*}$-algebra is a complete Hausdorff complex topological $*$\ -algebra $A$%
\ whose topology is determined by its continuous $C^{*}$-seminorms in the
sense that the net $\{a_i\}_{i\in I}$\ converges to $0$\ if and only if the
net $\{p(a_i)\}_{i\in I}$\ converges to $0$ for every continuous $C^{*}$%
-seminorm $p$\ on $A$. The terminology ''locally $C^{*}$-algebra'' is due to
Inoue (see [2]). Locally $C^{*}$-algebras were also studied by Phillips (
under the name of pro -$C^{*}$-algebra, see [7]), Fragoulopoulou, and other
people.

A representation of $A$ on a Hilbert space $H$ is a continuous $*$ -morphism 
$\varphi $ from $A$ to $L(H)$, the $C^{*}$-algebra of all bounded linear
operators on $H$. Given a locally $C^{*}$-algebra $A$ which acts
non-degenerately on a Hilbert module $E$ over a locally $C^{*}$-algebra $B$
and a non-degenerate representation $\left( \varphi ,H\right) $ of $B$,
exactly as in the case of $C^{*}$-algebras (see [8]), we construct a
representation of $A$, called the Rieffel-induced representation from $B$ to 
$A$ via $E$, and then we prove some properties of this representation. Thus,
we prove that the theorem on induction in stages (Theorem 5.9 in [8]) is
also true in the context of locally $C^{*}$-algebras (Theorem 3.6). In
section 4, we prove that if $A$ and $B$ are two locally $C^{*}$-algebras
which are strong Morita equivalent, then any non-degenerate representation
of $A$ is induced from a non-degenerate representation of $B$ (Theorem 4.4).

\section{Preliminaries}

Let $A$ be a locally $C^{*}$-algebra and let $S(A)$ be the set of all
continuous $C^{*}$-seminorms on $A$. If $p\in S(A)$, then $A_p=A/\ker p$ is
a $C^{*}$-algebra in the norm induced by $p$ and $A=\lim\limits_{\stackunder{%
p\in S(A)}{\leftarrow }}A_p$. The canonical map from $A$ onto $A_p$ is
denoted by $\pi _p$ and the image of $a$ under $\pi _p$ by $a_p$.

An isomorphism from a locally $C^{*}$-algebra $A$ to a locally $C^{*}$%
-algebra $B$ is a bijective, continuous $*$ -morphism $\Phi $ from A to $B$
such that $\Phi ^{-1}$ is continuous.

If $\left( \varphi ,H\right) $ is a representation of $A$, then there is $%
p\in S(A)$ and a representation $\left( \varphi _p,H\right) $ of $A_p$ such
that $\varphi =\varphi _p\circ \pi _p$. We say that $\left( \varphi
_p,H\right) $ is a representation of $A_p$ associated to $\left( \varphi
,H\right) $. The representation $\left( \varphi ,H\right) $ is
non-degenerate if $\varphi (A)H$ is dense in $H$. Clearly, $\left( \varphi
,H\right) $ is non-degenerate if and only if $\left( \varphi _p,H\right) $
is non-degenerate. We say that the representations $\left( \varphi
_1,H_1\right) $ and $\left( \varphi _2,H_2\right) $ of $A$ are unitarily
equivalent if there is a unitary operator $U$ from $H_1$ onto $H_2$ such
that $U\circ \varphi _1(a)$ $=\varphi _2(a)\circ U$ for all $a\in A$.

\begin{Definition}
A pre-Hilbert $A$-module is a complex vector space\ $E$ which is also a
right $A$-module, compatible with the complex algebra structure, equipped
with an $A$-valued inner product $\left\langle \cdot ,\cdot \right\rangle
:E\times E\rightarrow A\;$which is $\Bbb{C}$- and $A$-linear in its second
variable and satisfies the following relations: \ \ 

\ \ \ \ \ $\;\;(i)\;\;\left\langle \xi ,\eta \right\rangle ^{*}=\left\langle
\eta ,\xi \right\rangle \;\;$for every $\xi ,\eta \in E;$ \ \ \ \ \ \ \ \ 

\ \ \ \ \ \ $\;(ii)\;\left\langle \xi ,\xi \right\rangle \geq 0\;\;$ for
every $\xi \in E;$

\ \ \ \ \ \ \ $(iii)\;\left\langle \xi ,\xi \right\rangle =0\;\;$if and only
if $\xi =0.$ \ 

We say that $E\;$is a Hilbert $A$-module if $E\;$is complete with respect to
the topology determined by the family of seminorms $\{\left\| \cdot \right\|
_p\}_{p\in S(A)}$,$\;$where $\left\| \xi \right\| _p=\sqrt{p\left(
\left\langle \xi ,\xi \right\rangle \right) },\xi \in E\;\;$(Definition 4.1
of [7]).
\end{Definition}

Let $E\;$be a Hilbert $A$-module.\ For $p\in S(A)$,$\;$the vector space $%
E_p=E/\mathcal{E}_p$, where $\mathcal{E}_p=\{\xi \in E;p(\left\langle \xi
,\xi \right\rangle )=0\}$,$\;$is a Hilbert $A_p$-module\ with the action of $%
A_p$ on $E_p$ defined by $\left( \xi +\mathcal{E}_p\right) \left( a+\ker
p\right) =\xi a+\mathcal{E}_p$ and the inner product defined by $%
\left\langle \xi +\mathcal{E}_p,\eta +\mathcal{E}_p\right\rangle =\pi
_p\left( \left\langle \xi ,\eta \right\rangle \right) $ (Lemma 4.5 [7]).\
The canonical map from $E\;$onto $E_p$ is denoted by $\sigma _p$ and the
image of $\xi $ under $\sigma _p$ by $\xi _p$. Thus, for $p,q\in
S(A),\;p\geq q$, there is a canonical morphism of vector spaces $\sigma
_{pq}\;$from $E_p\;$ into $E_q\;$such that $\sigma _{pq}\left( \xi _p\right)
=\xi _q,\;\xi _p\in E_p$.$\;$Then $\{E_p,A_{p,}\sigma _{pq}:E_p\rightarrow
E_q,p\geq q;p,q\in S(A)\}$ is an inverse system of Hilbert $C^{*}$-modules
in the following sense: $\sigma _{pq}(\xi _pa_p)=\sigma _{pq}(\xi _p)\pi
_{pq}(a_p),\xi _p\in E_p,a_p\in A_p;$ $\left\langle \sigma _{pq}(\xi
_p),\sigma _{pq}(\eta _p)\right\rangle =\pi _{pq}(\left\langle \xi _p,\eta
_p\right\rangle ),\xi _p,\eta _p\in E_p;$ $\sigma _{pp}(\xi _p)=\xi _p,\;\xi
_p\in E_p\;$and $\sigma _{qr}\circ \sigma _{pq}=\sigma _{pr}\;$if $p\geq
q\geq r,$ and $\lim\limits_{\stackunder{p}{\leftarrow }}E_p$ is a Hilbert $A$%
-module which may be identified with the Hilbert $A$-module $E$ (Proposition
4.4 [7]).

A Hilbert $A$ -module $E$ is full if the ideal of $A$ generated by $%
\{\left\langle \xi ,\eta \right\rangle ,\xi ,\eta \in E\}$ is dense in $A$.

Let $E\;$and $F$\ be two Hilbert $A$-modules. The set of all adjointable
linear operators from $E$ to $F$ is denoted by $L_A(E,F)$, and we write $%
L_A(E)\;$for $L_A(E,E)$. We consider on $L_A(E,F)$ the topology determined
by the family of seminorms $\left\{ \widetilde{p}\right\} _{p\in S(A)}$ ,
where $\widetilde{p}(T)=\sup \left\{ \left\| T\xi \right\| _p;\left\| \xi
\right\| _p\leq 1\right\} $. Then $L_A(E,F)$ is isomorphic to $\lim\limits_{%
\stackunder{p}{\leftarrow }}L_{A_p}(E_p,F_p)$ (Proposition 4.7, [7]), and $%
L_A(E)$ becomes a locally $C^{*}$-algebra. The canonical maps from $L_A(E,F)$
to $L_{A_p}(E_p,$ $F_p)$, $p\in S(A)$ are denoted by $\left( \pi _p\right)
_{*}$ and $\left( \pi _p\right) _{*}\left( T\right) \left( \sigma _p(\xi
)\right) =\sigma _p\left( T\xi \right) $.

We say that the Hilbert $A$ -modules $E$ and $F$ are unitarily equivalent if
there is a unitary operator in $L_A(E,F)$.

A locally $C^{*}$-algebra $A$ acts non-degenerately on a Hilbert $B$ -module 
$E$ if there is a continuous $*$ -morphism $\Phi $ from $A$ to $L_B(E)$ such
that $\Phi (A)E$ is dense in $E$.

The closed vector subspace of $L_A(E,F)$ spanned by $\left\{ \theta _{\eta
,\xi };\xi \in E,\eta \in F\right\} $, where $\theta _{\eta ,\xi }(\zeta
)=\eta \left\langle \xi ,\zeta \right\rangle $, is denoted by $K_A(E,F)$,
and we write $K_A(E)$ for $K_A(E,E)$. Moreover, the locally $C^{*}$-algebras 
$K_A(E,F)$ and $\lim\limits_{\stackunder{p}{\leftarrow }}K_{A_p}(E_p,F_p)$
are isomorphic as well as the $C^{*}$-algebras $\left( K_A(E,F)\right) _p$
and $K_{A_p}(E_p,F_p)$ for all $p\in S(A)$. Since $K_A(E)E$ is dense in $E,$ 
$K_A(E)$ acts non-degenerately on $E$.

\section{Induced representations}

Let $A$ and $B$ be two locally $C^{*}$-algebras, let $E$ be a Hilbert $B$
-module, let $\Phi :A\rightarrow L_B(E)$ be a non-degenerate continuous $*$%
-morphism and let $\left( \varphi ,H\right) $ be a non-degenerate
representation of $B$. We will construct a non-degenerate representation $%
\left( _E^A\varphi ,_EH\right) $ of $A$ from $\left( \varphi ,H\right) $ via 
$E$.

\begin{Construction}
(for $C^{*}$-algebras, see [8]): Define a sesquilinear form $\left\langle
\cdot ,\cdot \right\rangle _0^{\varphi ^{}}$ on the vector space $E\otimes _{%
\text{alg}}H$ by 
\[
\left\langle \xi \otimes h_1,\eta \otimes h_2\right\rangle _0^{\varphi
^{}}=\left\langle h_1,\varphi \left( \left\langle \xi ,\eta \right\rangle
_E\right) h_2\right\rangle _{\varphi _{}} 
\]
where $\left\langle \cdot ,\cdot \right\rangle _{\varphi ^{}}$ denotes the
inner product on the Hilbert space $H$. It is easy to see that $\left(
E\otimes _{\text{alg}}H\right) /N_{\varphi ^{}}$, where $N_{\varphi ^{}}$ is
the vector subspace of $E\otimes _{\text{alg}}H$ generated by $\{\xi \otimes
h\in $ $E\otimes _{\text{alg}}H;\left\langle \xi \otimes h,\xi \otimes
h\right\rangle _0^{\varphi ^{}}=0\}$, is a pre-Hilbert space with the inner
product defined by 
\[
\left\langle \xi \otimes h_1+N_{\varphi ^{}},\eta \otimes h_2+N_{\varphi
^{}}\right\rangle ^{\varphi ^{}}=\left\langle \xi \otimes h_1,\eta \otimes
h_2\right\rangle _0^{\varphi ^{}}. 
\]
The completion of $\left( E\otimes _{\text{alg}}H\right) /N_{\varphi ^{}}$
with respect to the inner product $\left\langle \cdot ,\cdot \right\rangle
^{\varphi ^{}}$ is denoted by $_EH$. Let $T\in L_B(E)$. Define a linear map $%
_E\varphi \left( T\right) $ from $E\otimes _{\text{alg}}H$ into $E\otimes _{%
\text{alg}}H$ by 
\[
_E\varphi \left( T\right) \left( \xi \otimes h\right) =T\xi \otimes h. 
\]
If $\left( \varphi _q,H\right) $ is a representation of $B_q$ associated to $%
\left( \varphi ,H\right) $, then we have

$\left\langle _E\varphi \left( T\right) \left( \xi \otimes h\right)
,_E\varphi \left( T\right) \left( \xi \otimes h\right) \right\rangle
_0^{\varphi ^{}}=\left\langle h,\varphi \left( \left\langle T\xi ,T\xi
\right\rangle _E\right) h\right\rangle _{\varphi ^{}}$

$\;\;\;\;\;\;\;\;\;\;\;\;\;\;\;\;\;\;\;\;\;\;\;\;\;=\left\langle h,\varphi
_q\left( \left\langle \left( \pi _q\right) _{*}(T)\sigma _q(\xi ),\left( \pi
_q\right) _{*}(T)\sigma _q(\xi )\right\rangle _{E_q}\right) h\right\rangle
_{\varphi ^{}}$

$\;\;\;\;\;\;\;\;\;\;\;\;\;\;\;\;\;\;\;\;\;\;\;\;\;\leq \widetilde{q}\left(
T\right) \left\langle h,\varphi _q\left( \left\langle \sigma _q(\xi ),\sigma
_q(\xi )\right\rangle _{E_q}\right) h\right\rangle _{\varphi ^{}}$

$\;\;\;\;\;\;\;\;\;\;\;\;\;\;\;\;\;\;\;\;\;\;\;\;\;=\widetilde{q}\left(
T\right) \left\langle h,\left( \varphi _q\circ \pi _q\right) \left(
\left\langle \xi ,\xi \right\rangle _E\right) h\right\rangle _{\varphi ^{}}$

$\;\;\;\;\;\;\;\;\;\;\;\;\;\;\;\;\;\;\;\;\;\;\;\;\;=\widetilde{q}\left(
T\right) \left\langle \xi \otimes h,\xi \otimes h\right\rangle _0^{\varphi
^{}}\;$ $\;\;\;\;\;\;\;\;\;\;\;\;\;\;\;\;\;\;\;\;\;\;\;\;\;\;\;\;\;\;\;\;\;%
\;\;\;\;\;\;\;\;\;\;\;\;\;\;\;\;\;\;\;\;$\ $\;\;\;\;\;\;\;\;\;\;\;\;\;\;\;\;%
\;\;\;\;\;\;\;\;\;\;\;\;\;\;\;$\ $\;\;\;\;\;\;\;\;\;\;\;$ \ \ \ \ \ \ \ \ \
\ \ \ \ \ \ \ \ \ \ \ \ \ \ \ \ $\;\;\;\;\;\;\;\;\;\;\;\;\;\;\;\;\;\;\;\;\;%
\;\;\;\;\;\;\;\;\;\;\;\;\;\;\;$\ \ \ \ \ \ \ \ \ \ \ \ \ \ \ \ \ \ \ \ \ \ \
\ \ \ \ \ \ \ \ \ \ \ \ \ \ \ \ \ \ \ \ \ \ \ \ \ \ \ \ \ \ \ \ \ \ \ \ \ \
\ \ \ \ for all $\xi \in E$ and $h\in H$. From this we conclude that $%
_E\varphi \left( T\right) $ may be extended to a bounded linear operator $%
_E\varphi \left( T\right) $ on $_EH$. In this way we have obtained a map $%
_E\varphi $ from $L_B(E)\;$to $L\left( _EH\right) $. It is easy to see that $%
\left( _E\varphi ,_EH\right) $ is a representation of $L_B(E)$ on $_EH$.
Moreover, $_E\varphi $ is non-degenerate. Then $_E\varphi \circ \Phi $ is a
non-degenerate representation of $A$ on $_EH$ and it is denoted by $%
_E^A\varphi $.
\end{Construction}

\begin{Definition}
The representation $\left( _E^A\varphi ,_EH\right) $ constructed above is
called the Rieffel-induced representation from $B$ to $A$ via $E.$
\end{Definition}

\begin{Remark}
\begin{enumerate}
\item  Let $\left( \varphi _1,H_1\right) $ and $\left( \varphi _2,H_2\right) 
$ be two non-degenerate representations of $B$. If $\left( \varphi
_1,H_1\right) $ and $\left( \varphi _2,H_2\right) $ are unitarily
equivalent, then $\left( _E^A\varphi _1,_EH_1\right) $ and $\left(
_E^A\varphi _2,_EH_2\right) $ are unitarily equivalent.

\item  Let $F$ be a Hilbert $B$ -module which is unitarily equivalent to $E$%
. If $U$ is a unitary element in $L_B(E,F)$ and $A$ acts on $F$ by $%
a\rightarrow U\circ \Phi (a)\circ U^{*}$, then the representations $\left(
_E^A\varphi ,_EH\right) $and $\left( _F^A\varphi ,_FH\right) $of $A$ are
unitarily equivalent.
\end{enumerate}
\end{Remark}

\proof%
$(1):$ If $U$ is a unitary operator from $H_1$ onto $H_2$, then it is not
hard to check that the linear operator $V$ from $E\otimes _{\text{alg}}H_1$
onto $E\otimes _{\text{alg}}H_2$ defined by $V(\xi \otimes h)=\xi \otimes Uh$
may be extended to a unitary operator $V$ from $_EH_1$ onto $_EH_2$ and
moreover, $V\circ _E^A\varphi _1(a)=$ $_E^A\varphi _2(a)\circ V$ for all $a$
in $A$.

$(2):$ Consider the linear operator $W$ from $E\otimes _{\text{alg}}H$ onto $%
F\otimes _{\text{alg}}H$ defined by $W(\xi \otimes h)=U\xi \otimes h$. Then
we have 
\begin{eqnarray*}
\left( _F^A\varphi (a)\circ W\right) \left( \xi \otimes h\right) &=&\left(
U\circ \Phi (a)\circ U^{*}\right) \left( U\xi \right) \otimes h=U\left( \Phi
\left( a\right) \xi \right) \otimes h \\
&=&W\left( \Phi \left( a\right) \xi \otimes h\right) =\left( W\circ
_E^A\varphi (a)\right) \left( \xi \otimes h\right)
\end{eqnarray*}
for all $a$ in $A$, $\xi $ in $E$ and $h$ in $H$. It is not difficult to see
that $W$ may be extended to a unitary operator from $_EH$ onto $_FH$ and $%
_F^A\varphi (a)\circ W=W\circ _E^A\varphi (a)$ for all $a$ in $A$.%
\endproof%

\begin{Proposition}
Let $\left( \varphi ,H\right) $ be a non-degenerate representation of $B$.
If $\left( \varphi _q,H\right) $ is a non-degenerate representation of $B_q$
associated to $\left( \varphi ,H\right) $, then there is $p\in S(A)$ such
that $A_p$ acts non-degenerately on $E_q$ and the representations $\left(
_E^A\varphi ,_EH\right) $ and $\left( _{E_q}^{A_p}\varphi _q\circ \pi
_p,_{E_q}H\right) $ of $A$ are unitarily equivalent.
\end{Proposition}

\proof%
Define a linear map $U$ from $E\otimes _{\text{alg}}H$ into $E_q\otimes _{%
\text{alg}}H$ by 
\[
U\left( \xi \otimes h\right) =\sigma _q\left( \xi \right) \otimes h. 
\]
Since 
\begin{eqnarray*}
\left\langle U\left( \xi \otimes h\right) ,U\left( \xi \otimes h\right)
\right\rangle _0^{\varphi _q} &=&\left\langle h,\varphi _q\left(
\left\langle \sigma _q\left( \xi \right) ,\sigma _q\left( \xi \right)
\right\rangle _{E_q}\right) h\right\rangle _{\varphi _{}} \\
&=&\left\langle h,\left( \varphi _q\circ \pi _q\right) \left( \left\langle
\xi ,\xi \right\rangle _E\right) h\right\rangle _{\varphi _{}} \\
&=&\left\langle \xi \otimes h,\xi \otimes h\right\rangle _0^{\varphi ^{}}
\end{eqnarray*}
for all $\xi \in E$ and $h\in H$, $U$ may be extended to a bounded linear
operator $U$ from $_EH$ onto $_{E_q}H$. It is easy to verify that $U$ is
unitary and $U\circ _E\varphi (T)=\left( _{E_q}\varphi \circ \left( \pi
_q\right) _{*}\right) \left( T\right) \circ U$ for all $T\in L_B(E)$. Hence
the representations $\left( _E\varphi ,_EH\right) $ and $\left(
_{E_q}\varphi _q\circ \left( \pi _q\right) _{*},_{E_q}H\right) $ of $L_B(E)$
are unitarily equivalent.

The continuity of $\Phi $ implies that there is $p\in S(A)$ such that $%
\widetilde{q}(\Phi (a))\leq p(a)$ for all $a$ in $A$ and so there is a $*$
-morphism $\Phi _p$ from $A_p$ to $L_{B_q}(E_q)$ such that $\Phi _p\circ \pi
_p=\left( \pi _q\right) _{*}\circ \Phi $. Moreover, $\Phi _p$ is
non-degenerate. From

\begin{eqnarray*}
U\circ _E^A\varphi \left( a\right) &=&U\circ _E\varphi (\Phi (a))=\left(
_{E_q}\varphi _q\circ \left( \pi _q\right) _{*}\right) \left( \Phi
(a)\right) \circ U \\
&=&\left( _{E_q}\varphi _q\left( \Phi _p\left( \pi _p\left( a\right) \right)
\right) \right) \circ U=\left( _{E_q}^{A_p}\varphi _q\circ \pi _p\right)
\left( a\right) \circ U
\end{eqnarray*}
for all $a\in A$, we conclude that the representations $\left( _E^A\varphi
,_EH\right) $ and $(_{E_q}^{A_p}\varphi _q\circ \pi _p,$ $_{E_q}H)$ of $A$
are unitarily equivalent and the proposition is proved. 
\endproof%

\begin{Corollary}
If $\left( \varphi ,H\right) =\left( \tbigoplus\limits_{i\in I}\varphi
_i,\tbigoplus\limits_{i\in I}H_i\right) $, then $\left( _E^A\varphi
,_EH\right) $ is unitarily equivalent to $(\tbigoplus\limits_{i\in I}$ $%
_E^A\varphi _i,$ $\tbigoplus\limits_{i\in I}$ $_EH_i)$.
\end{Corollary}

\proof%
Let $\left( \varphi _q,H\right) $ be a representation of $B_q$ associated to 
$\left( \varphi ,H\right) $. It is easy to see that there is a
representation $\left( \varphi _{iq},H_i\right) $ of $B_q$ such that $%
\varphi _{iq}\circ \pi _q=\varphi _i$ for each $i\in I$. Moreover, $\varphi
_q=\tbigoplus\limits_{i\in I}\varphi _{iq}$. By Proposition 3.4, there is $%
p\in S(A)\;$such that the representations $\left( _E^A\varphi ,_EH\right) $
and $\left( _{E_q}^{A_p}\varphi _q\circ \pi _p,_{E_q}H\right) $ of $A$ are
unitarily equivalent as well as the representations $\left( _E^A\varphi
_i,_EH\right) $ and $\left( _{E_q}^{A_p}\varphi _{iq}\circ \pi
_p,_{E_q}H_i\right) $ for all $i\in I$.

On the other hand, we know that the representations $\left(
_{E_q}^{A_p}\varphi _q,_{E_q}H\right) $ and $\left( \tbigoplus\limits_{i\in
I}\text{ }_{E_q}^{A_p}\varphi _{iq},\tbigoplus\limits_{i\in I}\text{ }%
_{E_q}H_i\right) $ of $A_p$ are unitarily equivalent (Corollary 5.4 in [8]).
This implies that the representations $\left( _{E_q}^{A_p}\varphi _q\circ
\pi _p,_{E_q}H\right) $ and $\tbigoplus\limits_{i\in I}$ $%
(_{E_q}^{A_p}\varphi _{iq}\circ \pi _p,$ $\tbigoplus\limits_{i\in I}$ $%
_{E_q}H_i)$ of $A$ are unitarily equivalent and the corollary is proved. 
\endproof%

Let $A,$ $B$ and $C$ be three locally $C^{*}$-algebras, let $E$ be a Hilbert 
$B$ -module and $F$ a Hilbert $C$-module and let $\Phi _1:A\rightarrow
L_B(E) $ and $\Phi _2:B\rightarrow L_C(F)$ be non-degenerate continuous $*$-
morphisms. If $E\otimes _{\Phi _2}F$ is the inner tensor product of $E$ and $%
F$ using $\Phi _2$, then $E\otimes _{\Phi _2}F=\lim\limits_{\stackunder{r\in
S(C)}{\leftarrow }}E\otimes _{\Phi _{2r}}F_r$ and the locally $C^{*}$
-algebras $L_C(E\otimes _{\Phi _2}F)$ and $\lim\limits_{\stackunder{r\in S(C)%
}{\leftarrow }}L_{C_r}(E\otimes _{\Phi _{2r}}F_r)$ are isomorphic as well as 
$K_C(E\otimes _{\Phi _2}F)$ and $\lim\limits_{\stackunder{r\in S(C)}{%
\leftarrow }}K_{C_r}(E\otimes _{\Phi _{2r}}F_r)$, where $\Phi _{2r}=\left(
\pi _r\right) _{*}\circ \Phi _2$ (see [3]). Moreover, there is a
non-degenerate continuous $*$ -morphism $\left( \Phi _2\right) _{*}$ from $%
L_B(E)$ to $L_C(E\otimes _{\Phi _2}F)$ defined by $\left( \Phi _2\right)
_{*}(T)\left( \xi \otimes _{\Phi _2}\eta \right) =T\xi \otimes _{\Phi
_2}\eta $. Let $\Phi =\left( \Phi _2\right) _{*}\circ \Phi _1$. Then $\Phi $
is a non-degenerate continuous $*$ -morphism from $A$ to $L_C(E\otimes
_{\Phi _2}F)$.

\begin{Theorem}
Let $A,$ $B,C,E,F,$ $\Phi _1$ and $\Phi _2$ be as above. If $\left( \varphi
,H\right) $ is a non-degenerate representation of $C$, then the
representations $\left( _G^A\varphi ,_GH\right) $, where $G=E\otimes _{\Phi
_2}F$, and $\left( _E^A\left( _F^B\varphi \right) ,_E\left( _FH\right)
\right) $ of $A$ are unitarily equivalent.
\end{Theorem}

\proof%
Let $\left( \varphi _r,H\right) $ be a non-degenerate representation of $C_r$
associated to $\left( \varphi ,H\right) $. Then there is $q\in S(B)\;$and a
non-degenerate continuous $*$ -morphism $\Psi _{2q}:B_q\rightarrow
L_{C_r}(F_r)$ such that $\Psi _{2q}\circ \pi _q=\left( \pi _r\right)
_{*}\circ \Phi _2$ and there is $p\in S(A)\;$and a non-degenerate continuous 
$*$ -morphism $\Psi _{1p}:A_p\rightarrow L_{B_q}(E_q)$ such that $\Psi
_{1p}\circ \pi _p=\left( \pi _q\right) _{*}\circ \Phi _1$ and a
non-degenerate continuous $*$ -morphism $\Phi _p:A_p\rightarrow L_{C_r}(G_r)$
such that $\Phi _p\circ \pi _p=$ $\left( \pi _r\right) _{*}\circ \Phi $.

According to Proposition 3.4, the representations $\left( _G^A\varphi
,_GH\right) $ and $(_{G_r}^{A_p}\varphi _r\circ \pi _p,$ $_{G_r}H)$ of $A$
are unitarily equivalent as well as the representations $\left( _F^B\varphi
,_FH\right) $ and $(_{F_r}^{B_q}\varphi _r\circ \pi _q,$ $_{F_r}H)$ of $B$.
Since the representations $\left( _F^B\varphi ,_FH\right) $ and $%
(_{F_r}^{B_q}\varphi _r\circ \pi _q,_{F_r}H)$ of $B$ are unitarilly
equivalent, by Proposition 3.4 and Remark 3.3 (1) we deduce that the
representations $\left( _E^A\left( _F^B\varphi \right) ,_E\left( _FH\right)
\right) $ and $(_{E_q}^{A_p}\left( _{F_r}^{B_q}\varphi _r\right) \circ \pi
_p,$ $_{E_q}\left( _{F_r}H\right) )$ of $A$ are unitarily equivalent.

To show that the representations $\left( _G^A\varphi ,_GH\right) $ and $%
\left( _E^A\left( _F^B\varphi \right) ,_E\left( _FH\right) \right) $ of $A$
are unitarily equivalent it is sufficient to prove that the representations $%
\left( _{G_r}^{A_p}\varphi _r,_{G_r}H\right) $ and $\left(
_{E_q}^{A_p}\left( _{F_r}^{B_q}\varphi _r\right) ,_{E_q}\left(
_{F_r}H\right) \right) $ of $A_p$ are unitarily equivalent. But we know that
the representations $\left( _{X_r}^{A_p}\varphi _r,_{X_r}H\right) $, where $%
X_r=$ $E_q\otimes _{\Psi _{2q}}F_r$, and $(_{E_q}^{A_p}\left(
_{F_r}^{B_q}\varphi _r\right) ,$ $_{E_q}\left( _{F_r}H\right) )$ of $A_p$
are unitarily equivalent (Theorem 5.9 in [8]) and so it is sufficient to
prove that the representations $\left( _{X_r}^{A_p}\varphi _r,_{X_r}H\right) 
$ and $\left( _{G_r}^{A_p}\varphi _r,_{G_r}H\right) $ of $A_p$ are unitarily
equivalent.

It is not hard to check that the linear map $U:G_r\rightarrow X_r$ defined
by $U\left( \xi \otimes _{\Phi _{2r}}\eta \right) =\sigma _q\left( \xi
\right) \otimes _{\Psi _{2q}}\eta $ is a unitary operator in $%
L_{C_r}(G_r,X_r)$ and moreover, $\left( \Phi _{2r}\right) _{*}(T)=U^{*}\circ
\left( \Psi _{2q}\right) _{*}\left( \left( \pi _q\right) _{*}(T)\right)
\circ U$ for all $T$ in $L_B(E)$ (see the proof of Proposition 4.4 in [3]).
Since 
\begin{eqnarray*}
\Phi _p(\pi _p(a)) &=&\left( \pi _r\right) _{*}\left( \left( \Phi _2\right)
_{*}\left( \Phi _1\left( a\right) \right) \right) =\left( \Phi _{2r}\right)
_{*}\left( \Phi _1(a)\right) \\
&=&U^{*}\circ \left( \Psi _{2q}\right) _{*}\left( \left( \pi _q\right)
_{*}\left( \Phi _1(a)\right) \right) \circ U \\
&=&U^{*}\circ \left( \left( \Psi _{2q}\right) _{*}\circ \Psi _{1p}\right)
\left( \pi _p(a)\right) \circ U
\end{eqnarray*}
for all $a$ in $A$ and by Remark 3.3 (2), the representations $\left(
_{G_r}^{A_p}\varphi _r,_{G_r}H\right) $ and $\left( _{X_r}^{A_p}\varphi
_r,_{X_r}H\right) $ of $A_p$ are unitarily equivalent and the theorem is
proved.%
\endproof%

\section{The imprimitivity theorem}

Let $A$ and $B$ be locally $C^{*}$-algebras. We recall that $A$ and $B$ are
strongly Morita equivalent, written $A\thicksim _MB$, if there is a full
Hilbert $A$ -module $E$ such that the locally $C^{*}$ -algebras $B$ and $%
K_A(E)$ are isomorphic. The strong Morita equivalence is an equivalence
relation in the set of all locally $C^{*}$-algebras (see [4]). Also the
vector space $K_A(E,A)$, denoted by $\widetilde{E}$, is a full Hilbert $%
K_A(E)$ -module with the action of $K_A(E)$ on $K_A(E,A)$ defined by $\left(
T,S\right) \rightarrow T\circ S$ , $S\in K_A(E)$ and $T\in K_A(E,A)$, and
the inner product defined by $\left\langle T,S\right\rangle =T^{*}\circ S,$ $%
T,S\in K_A(E,A)$. Moreover, the linear map $\alpha $ from $A$ to $%
K_{K_A(E)}\left( \widetilde{E}\right) $ defined by $\alpha (a)\left( \theta
_{b,\xi }\right) =\theta _{ab,\xi }$ is an isomorphism of locally $C^{*}$%
-algebras (see [4]). Since the locally $C^{*}$-algebras $B$ and $K_A(E)$ are
isomorphic, $\widetilde{E}$ may be regarded as a Hilbert $B$ -module.

It is not hard to check that the linear operator $U_p$ from $\left( 
\widetilde{E}\right) _p$ to $\widetilde{E_p}$ defined by $U_p\left( T+\ker
\left( \widetilde{p}\right) \right) =\left( \pi _p\right) _{*}\left(
T\right) $ is unitary. Thus the Hilbert $K_{A_p}(E_p)$ -modules $\left( 
\widetilde{E}\right) _p$ and $\widetilde{E_p}$ may be identified.

\begin{Lemma}
If $A\thicksim _MB$, then for each $p\in S(A)$ there is $q_p\in S(B)\;$such
that $A_p\thicksim _MB_{q_p}.$ Moreover, the set $\{q_p\in S(B);p\in S(A)$
and $A_p\thicksim _MB_{q_p}\}$ is a cofinal subset of $S(B)$.
\end{Lemma}

\proof%
If $\Phi $ is an isomorphism of locally $C^{*}$-algebras from $B$ onto $%
K_A(E)$, then the map $\widetilde{p}\circ \Phi $, denoted by $q_p$, is a
continuous $C^{*}$-seminorm on $B$. Since $\ker \pi _{q_p}=\ker \left( \pi
_p\right) _{*}\circ \Phi $, there is a unique continuous $*$ -morphism $\Phi
_{q_p}$ from $B_{q_p}$ onto $K_{A_p}(E_p)$ such that $\Phi _{q_p}\circ \pi
_{q_p}=\left( \pi _p\right) _{*}\circ \Phi $. Moreover, $\Phi _{q_p}$ is an
isomorphism of $C^{*}$-algebras, and since $E_p$ is a full Hilbert $A_p$%
-module, we conclude that $A_p\thicksim _MB_{q_p}$.

To show that $\{q_p\in S(B);p\in S(A)$ and $A_p\thicksim _MB_{q_p}\}$ is a
cofinal subset of $S(B)$, let $q\in S(B)$. Then there is $p_0\in S(A)$ such
that 
\[
q\left( \Phi ^{-1}\left( \Phi \left( b\right) \right) \right) \leq 
\widetilde{p_0}\left( \Phi \left( b\right) \right) 
\]
for all $b\in B$, whence, since $q\left( \Phi ^{-1}\left( \Phi \left(
b\right) \right) \right) =q(b)$ and $\widetilde{p_0}\left( \Phi \left(
b\right) \right) =$ $q_{p_0}(b)$, we deduce that $q\leq q_{p_0}$.%
\endproof%

\begin{Remark}
If $E$ is a Hilbert $B$-module which gives the strong Morita equivalence
between the locally $C^{*}$-algebras $A$ and $B$, then $E_p$ gives the
strong Morita equivalence between the $C^{*}$-algebras $A_p$ and $B_{q_p}$.
\end{Remark}

\begin{Theorem}
Let $A$ and $B$ be two locally $C^{*}$-algebras such that $A\thicksim _MB$
and let $\left( \varphi ,H\right) $ be a non-degenerate representation of $A$%
. Then $\left( \varphi ,H\right) $ is unitarily equivalent to $\left( _{%
\widetilde{E}}^A\left( _E^B\varphi \right) ,_{\widetilde{E}}\left(
_EH\right) \right) $, where $E$ is a Hilbert $A$ -module which gives the
strong Morita equivalence between $A$ and $B$.
\end{Theorem}

\proof%
Let $\left( \varphi _p,H\right) $ be a non-degenerate representation of $A_p$
associated to $\left( \varphi ,H\right) $. By Lemma 4.1 there is $q\in S(B)$
such that $A_p\thicksim _MB_q$. Moreover, the Hilbert $A_p$ -module $E_p$
gives the strong Morita equivalence between $A_p$ and $B_q$ (Remark 4.2).
Then the representations $\left( \varphi _p,H\right) $ and $\left( _{%
\widetilde{E_p}}^{A_p}\left( _{E_p}^{B_q}\varphi _p\right) ,_{\widetilde{E_p}%
}\left( _{E_p}H\right) \right) $ of $A_p$ are unitarily equivalent (Theorem
6.23 in [8]) and by Remark 3.3 (2), the representations $\left( _{\widetilde{%
E_p}}^{A_p}\left( _{E_p}^{B_q}\varphi _p\right) ,_{\widetilde{E_p}}\left(
_{E_p}H\right) \right) $ and $\left( _{\widetilde{E}_p}^{A_p}\left(
_{E_p}^{B_q}\varphi _p\right) ,_{\widetilde{E}_p}\left( _{E_p}H\right)
\right) $ of $A_p$ are unitarily equivalent. From these facts we conclude
that the representations $\left( \varphi ,H\right) $ and $\left( _{%
\widetilde{E}_p}^{A_p}\left( _{E_p}^{B_q}\varphi _p\right) \circ \pi _p,_{%
\widetilde{E}_p}\left( _{E_p}H\right) \right) $ of $A$ are unitarily
equivalent.

On the other hand, according to Proposition 3.4, the representations $\left(
_E^B\varphi ,_EH\right) $ and $\left( _{E_p}^{B_q}\varphi _p\circ \pi
_q,_{E_p}H\right) $ of $B$ are unitarily equivalent. From this, using Remark
3.3(1) and Proposition 3.4, we deduce that the representations $\left( _{%
\widetilde{E}}^A\left( _E^B\varphi \right) ,_{\widetilde{E}}\left(
_EH\right) \right) $ and $\left( _{\widetilde{E}_p}^{A_p}\left(
_{E_p}^{B_q}\varphi _p\right) \circ \pi _p,_{\widetilde{E}_p}\left(
_{E_p}H\right) \right) $ of $A$ are unitarily equivalent and the theorem is
proved. 
\endproof%

\begin{Theorem}
Let $A$ and $B$ be locally $C^{*}$-algebras. If $A\thicksim _MB$, then there
is a bijective correspondence between equivalence classes of non-degenerate
representations of $A$ and $B$ which preserves direct sums and
irreducibility.
\end{Theorem}

\proof%
Let $E$ be a Hilbert $A$ -module which gives the strong Morita equivalence
between $A$ and $B$. By Theorem 4.3 and Remark 3.3 (1) the map from the set
of all non-degenerate representations of $A$ to the set of all
non-degenerate representations of $B$ which maps $\left( \varphi ,H\right) $
onto $\left( _E^B\varphi ,_EH\right) $ induces a bijective correspondence
between equivalence classes of non-degenerate representations of $A$ and $B$%
. Moreover, this correspondence preserves direct sums (Corollary 3.5).

Let $\left( \varphi ,H\right) $ be an irreducible, non-degenerate
representation of $A$. Suppose that $\left( _E^B\varphi ,_EH\right) $ is not
irreducible. Then $\left( _E^B\varphi ,_EH\right) =\left( \psi _1\oplus \psi
_2,H_1\oplus H_2\right) $ and by Corollary 3.5 and \ Theorem 4.3\ the
representations $(_{\widetilde{E}}^A\psi _1\oplus _{\widetilde{E}}^A\psi _2,$%
\ $_{\widetilde{E}}H_1\oplus _{\widetilde{E}}H_2)$ and $\left( \varphi
,H\right) $ of $A$ are unitarily equivalent, a contradiction. So the
bijective correspondence defined above preserves irreducibility.%
\endproof%

Department of Mathematics, Faculty of Chemistry, University of Bucharest,
Bd. Regina Elisabeta nr.4-12, Bucharest, Romania

\smallskip mjoita@fmi.unibuc.ro

\end{document}